\documentclass[preprint]{elsarticle}

\makeatletter
\def\ps@pprintTitle{}
\makeatother

\usepackage[a4paper,margin=1in]{geometry}
\usepackage{amsmath,amssymb,amsthm}
\usepackage{graphicx}
\usepackage{hyperref}
\usepackage{bm}
\usepackage{physics}
\usepackage{mathtools}
\usepackage{enumitem}

\theoremstyle{plain}

\theoremstyle{remark}

\theoremstyle{definition}

\title{Cayley-transform analysis and numerical validation of \\
the convergent Born series for the Helmholtz equation}

\tnotetext[1]{Preprint submitted to the Journal of Computational Physics.}

\author{Morten Jakobsen}

\address{Department of Earth Science, University of Bergen}

\begin{document}

\begin{frontmatter}


\begin{abstract}
We develop an operator-theoretic framework for the Convergent Born Series (CBS)
method applied to the Lippmann--Schwinger equation for high-frequency Helmholtz
problems. In contrast to the Fourier-based analysis of Osnabrugge et al., our
approach expresses the preconditioned Lippmann--Schwinger iteration entirely in
terms of the resolvent of a self-adjoint background operator. This leads to a
unitary Cayley-transform representation of the CBS iteration operator, from
which we derive basis-independent bounds on its numerical range and a general
convergence criterion valid on arbitrary bounded domains and for complex-valued
wave numbers. Because the analysis does not rely on an explicit Green's
function in the Fourier domain, the Cayley-transform framework extends
naturally to a broader class of frequency-domain wave and diffusion equations
whose fundamental solutions are not available in closed form. We further
incorporate smoothly tapered complex-wavenumber absorbing layers that preserve
the self-adjoint structure of the reference operator and enhance the
contractivity of the iteration without modifying the differential operator. In
addition to this theoretical generalization, we present a detailed numerical
validation in which CBS solutions are benchmarked against PML-based
finite-difference wavefield simulations. These experiments demonstrate that the
operator-theoretic CBS formulation delivers accurate and stable results across a
broad range of contrasts and frequencies, thereby significantly extending the
applicability and theoretical foundation of the CBS method beyond previously
analyzed settings.
\end{abstract}

\end{frontmatter}

\clearpage

\section{Introduction}

The Helmholtz equation plays a central role in modeling time--harmonic wave phenomena in acoustics, elastodynamics and electromagnetics.  In many applications the wavenumber varies spatially, giving rise to variable--coefficient Helmholtz problems that capture heterogeneities in the underlying medium \cite{weglein2003inverse,jakobsen2016renormalized}. The same mathematical structure appears in nonrelativistic quantum mechanics, where the time--independent Schrödinger equation can be written in Helmholtz form \cite{Sheng1995WaveScattering}. Accurate solutions to such heterogeneous Helmholtz problems are essential in areas ranging from seismic and ultrasonic imaging to photonics and computational quantum mechanics \cite{JakobsenXiangVanDongen2023, Sheng1995WaveScattering}. Because these settings often involve large, strongly varying media at high frequency, achieving numerical robustness and scalability remains a formidable challenge.

A wide range of numerical methods has been developed for solving the Helmholtz
equation, including classical finite--difference and finite--element schemes
\cite{marfurt1984,babuska2000}, as well as high--frequency discretizations such as
spectral--element \cite{komatitsch1998,komatitsch1999} and discontinuous Galerkin
methods \cite{feng2009}. Despite these advances, all approaches based on
differential operators face a fundamental difficulty at high wavenumbers: the
Helmholtz operator loses coercivity \cite{jakobsen2025biotallard}, making accurate numerical approximation
increasingly difficult. In finite--element methods, this behavior is reflected
in the well-known \emph{pollution effect}, where numerical dispersion forces the
number of degrees of freedom per wavelength to grow far beyond what basic
approximation arguments would suggest
\cite{babuska2000,ihlenburg1995}. Consequently, maintaining accuracy requires
progressively finer spatial resolution as the frequency increases. Even high-order
discretizations can display substantial errors at practically relevant mesh
sizes unless they are combined with very fine meshes or additional stabilization
\cite{cessenat1998}. 

Finite--difference methods encounter analogous challenges
\cite{marfurt1984,ajo2005frequency,AbubakarHabashy2013ViscoAcousticVIE,Huang}, and
their practical performance at high frequencies often depends on sophisticated
preconditioning strategies, such as shifted--Laplacian or Born-based approaches
\cite{Sheikh2016Shifted,GanderZhang2019,Huang}. As the wavenumber increases, wave
propagation and interference occur over many spatial scales, which reduces the
effectiveness of standard iterative solvers. In practice, this frequently leads
to slow convergence, stagnation, or even divergence
\cite{Frelet}, thereby limiting the scalability of
differential-operator-based formulations for large or strongly heterogeneous
problems. Such scalability limitations are particularly restrictive in seismic
applications, where large numbers of forward and adjoint solves are required, as
for example in acoustic full-waveform inversion
\cite{ajo2005frequency, Huang}.

These limitations strongly motivate the use of volume integral equation (VIE) formulations for high–frequency Helmholtz problem \cite{pike2001scattering, coltonkress1983, AbubakarHabashy2013ViscoAcousticVIE, osnabrugge2016, Krueger2017,malovichko2018acoustic, jakobsen2020renormalized}. Whereas differential methods produce large, indefinite systems, the Lippmann–Schwinger equation yields a linear system operator of the form $(I+K)$, with $K$ compact on suitable function spaces \cite{kleinman1991iterative, coltonkress2013}. This structure eliminates the high–frequency indefiniteness of the differential operator and keeps the essential spectrum tightly clustered at unity \cite{coltonkress2013,stevenson1990}. Consequently, the resulting systems are significantly more amenable to iterative solution, often requiring only mild preconditioning \cite{epstein2018, eikrem2021, ying2015, liu2018}. Moreover, the computational domain is restricted to the region of heterogeneity, and the free–space Green's function automatically enforces the radiation condition \cite{coltonkress2013}, enabling high–order accuracy and FFT–based matrix–free implementations \cite{epstein2018,greengard1987, jakobsen2016renormalized, jakobsen2020homotopy}.

Beyond these numerical advantages, scattering-series formulations provide a
transparent physical interpretation: each term corresponds to an additional
scattering event \cite{newton1992, pike2001scattering}. The Born series has been a central tool in scattering theory
since the 1920s, and its convergence difficulties have been examined for nearly
a century, underscoring the fundamental importance of understanding when and why
such series converge \cite{coltonkress2013, osnabrugge2016}. Depending on the scale, such series admit composite-particle
or effective-medium interpretations grounded in classical scattering theory \cite{Weinberg1963, Sheng1995WaveScattering, jakobsen2020renormalized}.
These representations naturally connect to renormalized and multiscale perspectives, and a convergent scattering series forms a compelling foundation for constructing inverse-scattering series. This perspective underlies a wide range of developments in direct and inverse scattering, including the inverse scattering series \cite{weglein2003inverse}, renormalized Lippmann--Schwinger approaches \cite{kouri2003renormalization}, and related renormalized or homotopy-based Born methods in seismic modeling \cite{jakobsen2016renormalized, jakobsen2020renormalized, jakobsen2020homotopy, jakobsen2020convergent}. Recent work in gravitational-wave theory has also exploited convergent Born expansions far beyond the traditional weak-scattering regime \cite{caronhuot2025bornseries}.

A major recent development in this direction is the \emph{Convergent Born Series} (CBS) introduced by Osnabrugge, Horstmeyer, and Vellekoop \cite{osnabrugge2016}. Their method achieves unconditional convergence for arbitrarily large and strongly scattering media through a simple pointwise complex preconditioner combined with a damped reference Green's function. However, their convergence analysis relies fundamentally on Fourier-symbol arguments and thus presupposes translation invariance and unbounded (or periodically extended) domains. As a result, the theoretical justification does not extend directly to bounded geometries, arbitrary boundary conditions, or heterogeneous backgrounds where no such spectral representation is available.

The goal of this paper is to develop an \emph{operator-theoretic formulation} of the Convergent Born Series that overcomes these limitations. By expressing the preconditioned Lippmann--Schwinger iteration entirely in terms of the resolvent of a self-adjoint background operator \cite{grahamSpence2019,coltonkress2013}, we derive a unitary Cayley-transform representation of the CBS iteration operator \cite{trefethenEmbree2005}. This yields a basis-independent convergence proof that applies on arbitrary bounded domains and does not rely on translation invariance or Fourier analysis. A key consequence is that the theory holds for general complex-valued wave numbers that describe dissipative media \cite{moiolaSpence2014,eikrem2021}, so physical attenuation and numerically introduced loss enter the analysis automatically through the imaginary part of the scattering potential.

This perspective also clarifies the role of absorbing boundary layers (ABLs). Since dissipation is treated naturally through complex wave numbers in the scattering potential, any smoothly varying complex extension of the medium parameters leads to additional contractivity of the CBS iteration. Motivated by this observation, we introduce a simple and effective ABL based on a smoothly tapered complex wavenumber \cite{jiang2004smoothcal}. This ABL preserves the self-adjoint structure of the reference operator, integrates seamlessly with the operator-theoretic CBS formulation and matrix-free integral-equation solvers, and provides a practical alternative to traditional PMLs \cite{berenger1994pml}, which typically modify the differential operator itself \cite{alles2011perfectly}. The inclusion of an efficient ABL within the operator-theoretic convergence analysis of the CBS iteration represents a novel and practically important feature, since efficient absorbing layers or PML constructions for integral-equation solvers remain comparatively rare in the literature; see, for example, \cite{alles2011perfectly}.

Within this operator framework, we adopt the local preconditioner introduced by
Osnabrugge et al.\ \cite{osnabrugge2016}, and show that it can be incorporated
into an operator-theoretic formulation of the preconditioned Lippmann--Schwinger
equation. By expressing the resulting fixed-point map in terms of the resolvent
of the background operator, we derive a compact quadratic representation of the
iteration using the Cayley transform. This formulation yields a clean estimate
on the numerical range of the iteration operator and an associated
basis-independent convergence criterion. In this way, our analysis extends
existing CBS convergence results to a significantly broader class of domains and
media, while retaining algebraic consistency with the Osnabrugge
preconditioner.

The remainder of the paper is organized as follows. Section~2 reviews the
Lippmann--Schwinger formulation together with both the conventional and the
preconditioned Born series, providing the multiple-scattering framework on which
our analysis builds. Section~3 develops the operator-theoretic foundation of the
paper: we introduce the relevant functional-analytic setting, derive the
Cayley-transform representation of the CBS iteration, and obtain a
basis-independent numerical-range estimate that leads to a general convergence
criterion. Section~4 turns to implementation aspects for the scalar Helmholtz
equation, including the construction of a viscoacoustic absorbing boundary layer
and an FFT-accelerated matrix-free realization of the CBS method. Section~5
presents numerical experiments validating the theory and comparing CBS to a
PML-based finite-difference solver. Section~6 discusses the physical
interpretation of the CBS and its scope and limitations in the scalar Helmholtz
setting. Section~7 provides concluding remarks.
 Appendix A summarizes the finite-difference frequency-domain (FDFD) solver with
perfectly matched layers (PML) used for numerical benchmarking. Appendix~B presents a theoretical extension of the operator-theoretic
CBS framework to frequency-domain variable-density acoustic problems. 
This appendix is included solely to illustrate the generality of the
convergence analysis at the level of operator theory.

\section{Aspects of multiple scattering theory}

\subsection{Conventional volume integral equation and Born series}

We consider the inhomogeneous scalar Helmholtz equation
\cite{coltonkress2013,osnabrugge2016,jakobsen2020convergent}
\begin{equation}
\nabla^{2}\psi(\mathbf{r}) + k(\mathbf{r})^{2}\,\psi(\mathbf{r}) = -S(\mathbf{r}), 
\qquad \mathbf{r} \in \mathbb{R}^{d},\ d\in\{2,3\},
\label{eq:LS-Helmholtz}
\end{equation}
where the source $S$ is compactly supported and the medium is assumed to be
gain-free, $\Im\,k(\mathbf{r})^{2} \ge 0$.
Equation~\eqref{eq:LS-Helmholtz} serves as the starting point for both
integral- and differential-equation approaches to wave propagation.

A common difficulty is the strongly indefinite character of the Helmholtz
operator.
To alleviate this, many iterative solvers introduce a small uniform complex
shift in the background wavenumber, the so-called \emph{shifted Laplacian}
\cite{Sheikh2016Shifted,GanderZhang2019}.
This modification moves the poles of the Green's function into the lower
half-plane and improves the behavior of the resolvent.
The same shift is used in the Convergent Born Series (CBS) of Osnabrugge
et~al.\ \cite{osnabrugge2016}, where it is interpreted as relocating the
residues of the free-space Green's function.

Introducing dissipation in the homogeneous reference medium is advantageous
from a numerical perspective, but it should not alter the physical scattering
problem \cite{Huang}.
We therefore rewrite the squared wavenumber in terms of a homogeneous reference
and a scattering potential,
\begin{equation}
k(\mathbf{r})^{2} = k_{0}^{2} + i\varepsilon + V(\mathbf{r}),
\qquad k_{0} > 0,\ \varepsilon > 0,
\label{eq:k2-split}
\end{equation}
where the modified scattering potential is defined as
\begin{equation}
V(\mathbf{r}) := k(\mathbf{r})^{2} - k_{0}^{2} - i\varepsilon.
\label{eq:V-definition}
\end{equation}
As pointed out by Osnabrugge et al.\ \cite{osnabrugge2016}, the final term in
$V$ constitutes a gain term that compensates for the artificial dissipation
introduced in the homogeneous reference medium, thereby preserving the
original physical scattering problem.

Let $g_\varepsilon$ denote the outgoing Green’s function of the damped operator
$\nabla^{2} + k_{0}^{2} + i\varepsilon$, satisfying
\begin{equation}
(\nabla^{2} + k_{0}^{2} + i\varepsilon)\,g_\varepsilon(\mathbf{r}) = -\delta(\mathbf{r}),
\label{eq:greens-kernel}
\end{equation}
with Fourier transform \cite{osnabrugge2016}
\begin{equation}
\widehat{g_\varepsilon}(\mathbf{k})
   = \frac{1}{\|\mathbf{k}\|^{2} - k_{0}^{2} - i\varepsilon}.
\label{eq:green-fourier}
\end{equation}
 The parameter $\varepsilon>0$ shifts the resolvent singularity away from the real axis and enforces the outgoing radiation condition; in Fourier space this corresponds to a small imaginary displacement of the pole, analogous to the $i\varepsilon$ prescription used to regularize on-shell singularities of massless propagators, while leaving the high-$\|\mathbf{k}\|$ asymptotics unchanged \cite{MandlShawQFT}.

We define the corresponding volume integral operator $G_\varepsilon$ by
\begin{equation}
(G_\varepsilon f)(\mathbf{r})
   := \int_{\mathbb{R}^{d}} g_\varepsilon(\mathbf{r}-\mathbf{r}')\, f(\mathbf{r}')\, \mathrm{d}r',
\qquad d\in\{2,3\},
\label{eq:greens-operator}
\end{equation}
which represents the resolvent of the Helmholtz operator in convolution form
\cite{coltonkress2013, osnabrugge2016}.

Applying $G_\varepsilon$ to \eqref{eq:LS-Helmholtz} and using the above
decomposition yields the Lippmann--Schwinger equation
\begin{equation}
\psi = G_\varepsilon V\psi + G_\varepsilon S,
\end{equation}
a standard formulation in scattering theory
\cite{morsefeshbach1953, newton1992, coltonkress2013, osnabrugge2016, jakobsen2020renormalized}.

This representation defines the classical Born fixed-point iteration
\[
\psi^{(n+1)} = G_\varepsilon V\psi^{(n)} + G_\varepsilon S,
\qquad
\psi^{(0)} := G_\varepsilon S,
\label{LSequation}
\]
which converges whenever $\rho(G_\varepsilon V) < 1$, generating the Born
(Neumann) series \cite{morsefeshbach1953, newton1992, coltonkress1983, osnabrugge2016}
\[
\psi = \sum_{m=0}^{\infty} (G_\varepsilon V)^{m} G_\varepsilon S
\]
In strongly scattering or extended media one typically has
$\rho(G_\varepsilon V) \ge 1$, so the unmodified fixed-point map is not
contractive, and the Born series fails to converge
\cite{newton1992, coltonkress2013, osnabrugge2016, jakobsen2020renormalized}.
We note that the damping parameter $\varepsilon>0$ does not generally reduce
$\rho(G_\varepsilon V)$, since it also enters $V$ with opposite sign.
This motivates the preconditioned formulation that underlies the Convergent
Born Series (CBS) \cite{osnabrugge2016}.

\subsection{Preconditioned volume integral equation and convergent Born series}

Following Osnabrugge et al.~\cite{osnabrugge2016}, we introduce a simple pointwise
preconditioning strategy that yields a uniformly convergent fixed-point scheme
for the Lippmann--Schwinger equation.
Let $\gamma$ denote a bounded diagonal operator.
Following Osnabrugge et al. \cite{osnabrugge2016}, we operator on \eqref{LSequation} by $\gamma$ from the left, add $\psi$ to both sides
and rearrange to obtain a preconditioned Lippmann-Schwinger equation:
\begin{equation}
\psi = M\psi + \gamma G_\varepsilon S,
\label{preLS_equation}
\end{equation}
with iteration operator
\begin{equation}
M := \gamma G_\varepsilon V - \gamma + I .
\label{eq:preconditioned-M}
\end{equation}
The convergence of the associated fixed-point iteration is governed entirely
by the spectral radius $\rho(M)$.

The key mechanism underlying the Convergent Born Series (CBS) introduced by Osnabrugge et al. \cite{osnabrugge2016} is the introduction
of a complex reference medium with squared wavenumber $k_0^2+i\varepsilon$,
which creates a balance between the uniform background loss $+i\varepsilon$ and
the compensating gain term $-i\varepsilon$ embedded in the scattering potential
$V$. Exploiting this structure, Osnabrugge et al.~\cite{osnabrugge2016} propose
the pointwise preconditioner
\begin{equation}
\gamma(\mathbf r) = \frac{i}{\varepsilon}\,V(\mathbf r).
\label{eq:gamma-choice}
\end{equation}
As discussed in the next section, substituting \eqref{eq:gamma-choice} into \eqref{eq:preconditioned-M} results in an iteration operator $M$ with very nice properties that will become very transparent when using a representation of the Green's function $G_{\epsilon }$ in terms of a unitary operator $U$.

A sufficient and practical admissibility condition ensuring contractivity is
\begin{equation}
\varepsilon \;\ge\; \|\Delta\|_\infty,
\qquad
\Delta(\mathbf r) := k(\mathbf r)^2 - k_0^2 ,
\label{eq:epsilon-admissibility}
\end{equation}
under which the spectral radius satisfies $\rho(M)<1$, provided that the medium
is gain-free, $\Im\,k(\mathbf r)^2 \ge 0$, and that absorption is present in at
least one region, either intrinsically or through an absorbing boundary layer.
Here $\|\cdot\|_\infty$ denotes the essential supremum over the computational
domain. Under these assumptions, condition \eqref{eq:epsilon-admissibility}
guarantees geometric convergence for arbitrarily large or strongly scattering
media and provides an explicit criterion for selecting $\varepsilon$ in
practice.

Under these conditions, the fixed-point scheme \eqref{preLS_equation} generates
the convergent Neumann expansion
\[
\psi = \sum_{m=0}^{\infty} M^{m}\bigl(\gamma G_\varepsilon S\bigr),
\]
which constitutes the Convergent Born Series (CBS), the preconditioned analogue
of the classical Born series, and remains convergent far beyond the
weak-scattering regime. As discussed by Osnabrugge et al. \cite{osnabrugge2016}, the scattering series solution of the Helmholtz equation shown above convergences for arbitrarily large media with arbitrarily strong scattering potentials. Thus, the scalar CBS is applicable to seismic wavefield modeling in strongly scattering media \cite{jakobsen2020renormalized} and other applications far beyond it's original domain of applicability, which was optics \cite{osnabrugge2016}.

The admissibility condition \eqref{eq:epsilon-admissibility} is justified and
extended in a basis-independent operator-theoretic framework in Section~3,
while practical choices of $\varepsilon$ and the role of absorbing boundary
layers are discussed in Section~4.

\section{Cayley Transform Analysis of the CBS Iteration}
\label{sec:operator_analysis}

\subsection{Functional‑analytic setting}
\label{sec:41}

Rather than restricting our analysis to the scalar Helmholtz equation (1), we here consider a more general class of scalar og vectorial linear PDE's written abstractly as
\begin{equation}
L\psi = -S,
\label{eq:PDE-general}
\end{equation}
where $L$ is a linear differential operator acting on a field $\psi$ and
$S\in\mathcal{H}$ is the source.

For the CBS iteration we assume that $L$ admits the decomposition
\begin{equation}
L = A + V,
\label{eq:L-decomposition}
\end{equation}
where $A$ is a \emph{self-adjoint} reference operator and $V$ is a
\emph{bounded pointwise multiplication operator} encoding all material
heterogeneity. Recall that a self-adjoint operator is a linear operator whose adjoint is equal to the operator itself and this ensure the spectrum is real. 

This structure is standard in the scalar Helmholtz equation,
$A=\nabla^{2}+k_{0}^{2}$ and $V(\mathbf{r})=k(\mathbf{r})^{2}-k_{0}^{2}$.  
More generally, many second‑order wave and diffusion equations
(heterogeneous acoustics, elasticity, Maxwell equations, poroelasticity,
anisotropic diffusion, etc.) can be written in this form after a symmetric
first‑order reformulation in which all spatial derivatives appear only
through fixed operators, and the heterogeneity appears entirely in a
multiplicative scattering potential \cite{deHoop1995Handbook, jakobsenENUMATHCoupledIE, JakobsenXiangVanDongen2023, saputeraEMdomainDecomposition, jakobsenTransitionOperatorFWI}.

Throughout we work in the Hilbert space
\[
\mathcal{H}=L^{2}(\Omega;\mathbb{C}^{m}),\qquad
\langle u,v\rangle = \int_{\Omega} u(\mathbf{r})^{*} v(\mathbf{r})\,dr,
\]
where $m=1$ corresponds to scalar acoustics and $m>1$ includes
multi-component systems such as vectorial acoustics
\cite{jakobsenENUMATHCoupledIE, JakobsenXiangVanDongen2023},
electromagnetics \cite{Krueger2017, saputeraEMdomainDecomposition},
elastodynamics \cite{jakobsenTransitionOperatorFWI},
and coupled diffusion–wave dynamics \cite{jakobsen2025biotallard}.
The operator $A$ is assumed self-adjoint under the imposed conservative
boundary conditions, while absorbing boundary layers and other dissipative
effects are represented through the scattering potential operator $V$.
In first-order form, this self-adjointness corresponds to energy conservation
in lossless media. It is also possible to use skew-adjoint operators $A$,
since multiplication by the imaginary unit $i$ renders them self-adjoint.
This formulation makes the analogy between classical wave theory and the
Schrödinger equation more transparent \cite{Sheng1995WaveScattering}.

For any $\varepsilon>0$ we define the shifted operator
\[
H_{\varepsilon} := A + i\varepsilon I.
\]
Since $\sigma(A)\subset\mathbb{R}$, the point $i\varepsilon$ lies outside the
spectrum, and consequently $H_{\varepsilon}$ is bijective, meaning that for
every right-hand side $f\in\mathcal H$ the equation $H_{\varepsilon}u=f$
admits a unique solution $u\in\mathcal H$ \cite{ReedSimon1972, coltonkress1983, coltonkress2013}.
The associated inverse operator
\[
G_{\varepsilon} := H_{\varepsilon}^{-1}
\]
is the shifted resolvent used in the CBS fixed-point map; the imaginary shift
introduces weak background damping, which guarantees uniqueness of the outgoing
field and suppresses resonant behavior \cite{newton1992, coltonkress2013}.

The contrast operator acts pointwise as
\[
(Vu)(\mathbf{r})=V(\mathbf{r})u(\mathbf{r}),\qquad V(\mathbf{r})\in \mathbb{C}^{m\times m},\quad
V\in L^{\infty}(\Omega),
\]
so $G_{\varepsilon}V$ is bounded.  For $m=1$ this gives a complex scalar
contrast; for $m>1$, $V(\mathbf{r})$ represents perturbations of material tensors
such as density, compressibility, permittivity, stiffness, or
conductivity.  No additional structure beyond boundedness is required. The handbook of scattering and radiation by de Hoop \cite{deHoop1995Handbook} shows how to represent different PDE's by equivalent integral equations of the Lippmann-Schwinger type required to use the operator theory present here. 
As discussed in the papers dealing with vectorial acoustics
\cite{jakobsenENUMATHCoupledIE, JakobsenXiangVanDongen2023},
electromagnetics \cite{Krueger2017, saputeraEMdomainDecomposition},
elastodynamics \cite{jakobsenTransitionOperatorFWI}, one can always ensure that the scattering potential is a pointwise multiplication operator if one formulates the integral equation in a higher-dimensional space. See appendix B for a detailed discussion of the operator structure in the variable-density acoustic case. 

These assumptions ensure that all operators in the preconditioned CBS
iteration are well defined and bounded on $\mathcal{H}$, providing the
functional‑analytic foundation for the Cayley‑transform analysis in the next subsection.

\subsection{Derivation of the CBS operator via the Cayley transform}
\label{sec:CBS-Cayley-derivation}

We define the Cayley transform of a self-adjoint operator $A$ by
\begin{equation}
U := (A - i\varepsilon I)(A + i\varepsilon I)^{-1},
\label{eq:Cayley-def-main}
\end{equation}
where $\varepsilon>0$. Since $A=A^*$, the operators $A\pm i\varepsilon I$ are boundedly invertible.
Taking adjoints and using $(BC)^*=C^*B^*$, $(B^{-1})^*=(B^*)^{-1}$, and $(A\pm i\varepsilon I)^*=A\mp i\varepsilon I$, we obtain
\[
U^*
= \bigl((A+i\varepsilon I)^{-1}\bigr)^*(A-i\varepsilon I)^*
= (A - i\varepsilon I)^{-1}(A + i\varepsilon I).
\]
A direct calculation then shows that $UU^*=U^*U=I$, and hence $U$ is unitary.

Spectrally, the Cayley transform maps each $\lambda\in\sigma(A)\subset\mathbb{R}$
to
\[
\frac{\lambda-i\varepsilon}{\lambda+i\varepsilon},
\]
which lies on the unit circle. Thus $U$ provides a unitary representation of the
background operator that captures its spectral geometry independently of the
domain, coordinate system, and boundary conditions.

From the definition of $U$ in equation ~\eqref{eq:Cayley-def-main} we compute
\begin{align*}
I-U
&= (A+i\varepsilon I)(A+i\varepsilon I)^{-1}
   - (A-i\varepsilon I)(A+i\varepsilon I)^{-1} \\
&= 2i\varepsilon\,(A+i\varepsilon I)^{-1}.
\end{align*}
Rearranging yields the exact identity
\begin{equation}
G_{\epsilon } = 
(A+i\varepsilon I)^{-1}
= \frac{1}{2i\varepsilon}(I-U).
\label{eq:resolvent-Cayley-main}
\end{equation}
This formula shows that the resolvent of the background operator depends on $A$
only through the single unitary operator $U$. This is a major advantages since unitary operators are much easier to deal with than arbitrary operators in convergence analysis. 

It follows from equations (9) and (10) that the CBS iteration operator has the form
\begin{equation}
M
= I - \frac{i}{\varepsilon} V
+ \frac{i}{\varepsilon} VG_{\epsilon } V,
\label{eq:CBS-def-main}
\end{equation}
 Substituting the expression for the Green's function $G_{\epsilon }$
from equation \eqref{eq:resolvent-Cayley-main} into the quadratic term in the above expression for $M$ gives
\begin{align*}
\frac{i}{\varepsilon} V G_{\epsilon } V
&= \frac{i}{\varepsilon}
   V \left( \frac{1}{2i\varepsilon}(I-U) \right) V \\
&= \frac{1}{2\varepsilon^{2}}
   \bigl( V^{2} - VUV \bigr).
\end{align*}
Inserting this expression for the quadratic term into the expression ~\eqref{eq:CBS-def-main} for the CBS iteration operator $M$ and collecting terms, we
obtain
\begin{equation}
M =
\frac{1}{2\varepsilon^{2}}
\bigl(
   -V^{2}
   + VUV
   - 2 i \varepsilon V
   + 2\varepsilon^{2} I
\bigr).
\label{eq:M-Cayley-main}
\end{equation}
The above equation is similar to equation (A.3) in  Osnabrugge \emph{et~al.}~\cite{osnabrugge2016}. 
In that work, the unitary operator $U$ is defined in a specific Fourier basis associated with the scalar Green’s function. 
In contrast, our expression~(13) for $U$ is formulated for a general self-adjoint operator~$A$, without reference to a particular basis or underlying scalar model. 
This abstraction is essential for our purposes, which are not merely to reproduce the convergence analysis of Osnabrugge \emph{et~al.}, but to extend it in a way that allows the coherent Born series to be applied to more general wave and diffusion systems.

Formula~\eqref{eq:M-Cayley-main} clearly separates the roles of the background and the contrast. 
All geometric and nonlocal effects enter exclusively through the unitary operator $U$, 
while $V$ remains a bounded multiplicative contrast; the only nonlocal contribution is $VUV$, encoding how the background couples the heterogeneous region to itself.

\subsection{Derivation of a global contractivity condition}

In order to prove convergence of the CBS iteration, it is sufficient to show that the spectral radius of the iteration operator $M$ satisfies $\rho(M)<1$. A sufficient condition for this is that the numerical range of $M$ is contained in the open unit disk. Following \cite{osnabrugge2016}, we therefore require that
\begin{equation}
\left|\langle \psi, M\psi\rangle\right|
< \|\psi\|^{2},
\qquad
\psi \in L^{2}(\Omega;\mathbb{C}^{m}),\ \psi \neq 0 .
\label{numRangeM}
\end{equation}
This numerical-range condition is sufficient, though not necessary, for $\rho(M)<1$, but has the advantage of yielding a robust, basis-independent convergence criterion for the CBS iteration. 
Although the intermediate estimates below are expressed using non-strict inequalities, the presence of absorption—either physical or through an absorbing boundary layer—rules out the marginal elastic case $\rho(M)=1$ identified by \cite{osnabrugge2016} and ensures strict contractivity of the CBS iteration.

Combining equations \eqref{numRangeM} and \eqref{eq:M-Cayley-main}, we get
\[
\langle \psi, M\psi\rangle
=
\frac{1}{2\varepsilon^{2}}
\Bigl(
\langle \psi, (-V^{2}-2 i \varepsilon V + 2\varepsilon^{2}I)\psi\rangle
+
\langle \psi, VUV\psi\rangle
\Bigr).
\]
Taking absolute values and applying the triangle inequality, we get
\begin{equation}
|\langle \psi, M\psi\rangle|
\le
\frac{1}{2\varepsilon^{2}}
\Bigl(
\bigl|
\langle \psi, (-V^{2}-2 i \varepsilon V + 2\varepsilon^{2}I)\psi\rangle
\bigr|
+
\bigl|
\langle \psi, VUV\psi\rangle
\bigr|
\Bigr).
\label{triangleEstimate}
\end{equation}
All terms in the quadratic form are local except for the mixed contribution
$\langle \psi, VUV\psi\rangle$.
This term can be controlled using elementary Hilbert-space estimates.
By the Cauchy--Schwarz inequality,
\[
\bigl|\langle \psi, VUV\psi\rangle\bigr|
\le
\|\psi\|\,\|VUV\psi\|.
\]
Since $U$ is unitary on $L^{2}(\Omega;\mathbb{C}^{m})$ and $V$ is a bounded
multiplication operator, the operator norm is submultiplicative, and hence
\begin{equation}
\|VUV\psi\|
\le
\|V\|\,\|UV\psi\|
=
\|V\|\,\|V\psi\|
\le
\|V\|^{2}\,\|\psi\|.
\label{VUV}
\end{equation}
Combining these inequalities gives
\begin{equation}
\bigl|\langle \psi, VUV\psi\rangle\bigr|
\le
\|V\|^{2}\,\|\psi\|^{2},
\label{eq:mixed-term-bound}
\end{equation}
without imposing any self-adjointness, normality, or positivity assumptions on $V$.
Substituting~\eqref{eq:mixed-term-bound} into the triangle-inequality
estimate \eqref{triangleEstimate} above therefore yields
\begin{equation}
\bigl| \langle \psi, M\psi\rangle \bigr|
\le
\frac{1}{2\varepsilon^{2}}
\Bigl(
\bigl|
\langle \psi,
(-V^{2}-2 i \varepsilon V + 2\varepsilon^{2}I)\psi
\rangle
\bigr|
+
\|V\|^{2}\,\|\psi\|^{2}
\Bigr),
\label{eq:NR-first-bound}
\end{equation}
which can be referred to as a global contractivity constraint.

\subsection{From global to local contractivity condition}
\label{subsec:pointwise-constraint}

Let us now try to reduce the global contractivity constraint
\eqref{eq:NR-first-bound} to a purely local condition on the scattering
potential. Since $V$ is a local multiplication operator in $L^{2}(\Omega;\mathbb{C}^{m})$, it follows that the scalar product term in equation \eqref{eq:NR-first-bound} can be expressed as 
\[
\bigl|
\langle \psi,
(-V^{2}-2 i \varepsilon V + 2\varepsilon^{2}I)\psi
\rangle
\bigr|
=
\left|
\int_{\Omega}
\psi(\mathbf{r})^{*}
\bigl(
- V(\mathbf{r})^{2}
- 2 i \varepsilon V(\mathbf{r})
+ 2\varepsilon^{2}I
\bigr)
\psi(\mathbf{r})\,d\mathbf{r}
\right|.
\]
Using the triangle inequality for integrals,
$
\left|
\int_{\Omega} f(\mathbf{r})\,d\mathbf{r}
\right|
\le
\int_{\Omega} |f(\mathbf{r})|\,d\mathbf{r},
$
followed by the elementary matrix bound
$|x^{*}Ax| \le \|A\|\,\|x\|^{2}$, we obtain
\begin{equation}
\bigl|
\langle \psi,
(-V^{2}-2 i \varepsilon V + 2\varepsilon^{2}I)\psi
\rangle
\bigr|
\le
\int_{\Omega}
\bigl\|
- V(\mathbf{r})^{2}
- 2 i \varepsilon V(\mathbf{r})
+ 2\varepsilon^{2}I
\bigr\|
\,\|\psi(\mathbf{r})\|^{2}\,d\mathbf{r}.
\label{AboveEstimate}
\end{equation}
Introducing the shifted contrast
\[
\Delta(\mathbf{r}) := V(\mathbf{r}) + i\varepsilon I,
\]
it follows from the global contractivity condition \eqref{eq:NR-first-bound} and the above estimate \eqref{AboveEstimate} that
\begin{equation}
\bigl|\langle \psi, M\psi\rangle\bigr|
\le
\frac{1}{2\varepsilon^{2}}
\int_{\Omega}
\Bigl(
\bigl\|
\varepsilon^{2}I - \Delta(\mathbf{r})^{2}
\bigr\|
+
\|V(\mathbf{r})\|^{2}
\Bigr)
\|\psi(\mathbf{r})\|^{2}\,d\mathbf{r}.
\label{GlobalIntegral}
\end{equation}
Since $\psi$ is arbitrary, it follows from the above version of the global contractivity condition \eqref{GlobalIntegral} that 
\begin{equation}
\bigl\|
\varepsilon^{2} I - \Delta(\mathbf{r})^{2}
\bigr\|
+
\|V(\mathbf{r})\|^{2}
\;\le\;
2\varepsilon^{2},
\qquad
\mathbf{r}\in\Omega.
\label{eq:pointwise-matrix-ineq}
\end{equation}
Condition \eqref{eq:pointwise-matrix-ineq} is a purely local matrix
inequality depending only on the pointwise values of the scattering
potential operator or matrix $V(\mathbf{r})$. All nonlocal effects arising from wave propagation have been
eliminated at this stage and enter the analysis only indirectly through
the unitary bound established in Section~3.3.

\subsection{From local to scalar contractivity condition}

To obtain a scalar admissibility criterion, we reduce the pointwise matrix
inequality \eqref{eq:pointwise-matrix-ineq} to a condition on the eigenvalues
of $V(\mathbf{r})$. Since an operator--norm bound implies the corresponding
quadratic--form inequality on all unit vectors, \eqref{eq:pointwise-matrix-ineq}
must hold in particular when evaluated on eigenvectors of $V(\mathbf{r})$.

Let $x$ be a unit eigenvector with eigenvalue $\lambda$, so that
$V(\mathbf{r})x=\lambda x$. Because $\Delta(\mathbf{r})=V(\mathbf{r})+i\varepsilon I$,
we have $\Delta(\mathbf{r})x=(\lambda+i\varepsilon)x$, and hence
\eqref{eq:pointwise-matrix-ineq} implies the scalar inequality
\begin{equation}
\bigl|
\varepsilon^{2} - (\lambda+i\varepsilon)^{2}
\bigr|
+|\lambda|^{2}
\;\le\;
2\varepsilon^{2}.
\label{eq:eigenvalue-ineq}
\end{equation}
A direct calculation shows that\eqref{eq:eigenvalue-ineq} holds if  
\begin{equation}
|\lambda(\mathbf{r})| \le \varepsilon,
\qquad
\Im \lambda(\mathbf{r}) \ge 0,
\qquad
\mathbf{r}\in\Omega.
\label{eq:admissibility-final}
\end{equation}
Condition \eqref{eq:admissibility-final} is therefore a sufficient pointwise
admissibility condition for \eqref{eq:pointwise-matrix-ineq}, and coincides
with the scalar criterion obtained by Osnabrugge et al. ~\cite{osnabrugge2016}.

\section{Some details for implementation}
\label{sec:ABL-implementation}

\subsection{Viscoacoustic Absorbing Boundary Layer}
\label{sec:ABL}

Artificial reflections arising from domain truncation present a major source of error in numerical simulations of wave propagation. To suppress these non-physical effects, we surround the physical domain with an absorbing boundary layer (ABL) that emulates a lossy exterior medium. The use of an effective ABL is especially critical in frequency-domain formulations, where the global nature of the solution precludes truncation of the wavefield prior to boundary interaction. Beyond its practical role in controlling reflections, the ABL must therefore be incorporated consistently into the convergence analysis. Since our formulation allows for a dissipative scattering potential operator \(V\), it is natural to include the ABL in \(V\) rather than in the reference operator \(A\), which is required to remain self-adjoint and non-dissipative. This separation preserves the structural properties of \(A\) while accounting for attenuation and boundary absorption entirely through \(V\), ensuring compatibility with the underlying operator framework and its convergence properties.

Our approach is based on a \emph{constant-$Q$ model}, in which attenuation is 
introduced by allowing the wave speed to become complex-valued \cite{kjartansson1979constantQ, carcione2022wave}.  We define the 
inverse quality factor as
\begin{equation}
Q^{-1}(\mathbf{r}) = Q^{-1}_{\max}\, T(\mathbf{r}),
\label{eq:Qinv_def}
\end{equation}
where $Q^{-1}_{\max}$ is a user-specified maximum attenuation strength and 
$T(\mathbf{r})$ is a smooth taper function taking values in $[0,1]$. The taper increases 
monotonically from zero at the interface between the physical domain and the ABL 
to unity at the outer boundary of the layer.

Any sufficiently smooth monotone taper is admissible; our convergence analysis 
places no restrictions on its specific functional form. For concreteness, and 
because it performs well in practice, we employ the standard cosine--squared 
taper
\begin{equation}
T(r)
  = \sin^{2}\!\left(
      \frac{\pi}{2}\,
      \frac{r - r_{\mathrm{in}}}{r_{\mathrm{out}} - r_{\mathrm{in}}}
    \right),
\qquad
r_{\mathrm{in}} \le r \le r_{\mathrm{out}},
\label{eq:taper}
\end{equation}
with $T(r)=0$ for $r \le r_{\mathrm{in}}$ and $T(r)=1$ for $r \ge r_{\mathrm{out}}$.
This choice is merely illustrative; polynomial tapers and other smooth profiles 
yield comparable performance.

The complex wave speed inside the ABL is defined by
\begin{equation}
C_{\mathrm{ABL}}(\mathbf{r})
  = C(\mathbf{r})\,\sqrt{1 - i\,Q^{-1}(\mathbf{r})},
\label{eq:CABL}
\end{equation}
where $C(\mathbf{r})$ is the real-valued wave speed in the absence of attenuation. 
This leads directly to the complex wavenumber
\begin{equation}
k(\mathbf{r}) = \frac{\omega}{C_{\mathrm{ABL}}(\mathbf{r})}.
\label{eq:kABL}
\end{equation}
This formulation is consistent with frequency-domain constant-$Q$ viscoacoustic 
models \cite{kjartansson1979constantQ} but does not rely on a viscoelastic 
modulus as in the classical Voigt model \cite{carcione2022wave}.  Instead, attenuation is introduced 
solely through the wave speed, making this approach especially well suited to scalar 
integral-equation solvers such as the Convergent Born Series (CBS). Because 
the attenuation enters exclusively through the multiplicative contrast in the 
Lippmann--Schwinger equation, the self-adjoint structure of the reference 
operator is preserved.

Unlike perfectly matched layers (PMLs) developed for differential formulations
\cite{berenger1994pml,chew1996pml,nataf2013pml}, which rely on complex coordinate
stretching and typically introduce auxiliary variables, the proposed ABL acts
solely through a modification of the physical medium. PML constructions tailored
to frequency-domain integral equations have also been introduced
\cite{alles2011perfectly}; however, these approaches remain based on coordinate
stretching and do not preserve the self-adjoint reference operator required by
the operator-theoretic CBS analysis. Moreover, the perfectly matched layer
proposed in \cite{alles2011perfectly} necessitates a variable-density
formulation. In contrast, the present ABL introduces attenuation in a simple and
direct manner that preserves the self-adjoint nature of the background operator
and integrates naturally with the Cayley-transform convergence framework
developed in Section~3.

\subsection{FFT-accelerated matrix-free implementation of the CBS}
\label{sec:implementation}

Our implementation follows the matrix-free Fourier-accelerated scheme \cite{osnabrugge2016, jakobsen2020renormalized, jakobsen2020homotopy}, but is presented here at the 
continuous operator level rather than in discrete notation.  
This formulation naturally accommodates a genuinely complex-valued wavenumber 
field and therefore allows the dissipative absorbing boundary layer (ABL) of 
Section~4.1 to be incorporated directly through the imaginary part of the 
scattering potential $V$, without modifying the background operator $H$ or the 
structure of the FFT-based resolvent.

A key practical consideration is the choice of the convergence-control parameter
$\varepsilon$. The analysis in Section~3 establishes that the CBS iteration is
contractive provided
\[
\varepsilon \;\ge\; \|\Delta\|_\infty, 
\qquad 
\Delta(\mathbf r) = k(\mathbf r)^2 - k_0^2,
\]
and we therefore adopt this admissibility condition directly. For each frequency,
the quantity $\|\Delta\|_\infty$ is evaluated once on the discrete computational
grid, yielding a value of $\varepsilon$ that naturally incorporates the effects
of both the physical medium and the absorbing boundary layer (ABL) introduced in
Section~4.1.

In practice, choosing $\varepsilon = 1.1\,\|\Delta\|_\infty$ provides a reliable
balance between damping and numerical stability. Smaller values of $\varepsilon$
lead to a sharper spectral kernel $\widehat g_\varepsilon$, while larger values
increase smoothing and strengthen contractivity. This strategy closely follows
the parameter-selection procedure of Osnabrugge et al.~\cite{osnabrugge2016};
however, in the present formulation $\Delta(\mathbf r)$ explicitly includes the
ABL, resulting in minor but important differences between our implementation and
that of~\cite{osnabrugge2016}.

The implementation is matrix-free and uses the spectral evaluation of convolution
with the outgoing damped Green's function.  
Let $\widehat g_\varepsilon(\mathbf k)$ denote the Fourier transform of this
Green’s function (see Section~2).
For any \emph{generic test function} $u$, the background resolvent is applied as
\[
G_\varepsilon u(\mathbf r)
=
\mathcal F^{-1}\!\bigl(
\widehat g_\varepsilon(\mathbf k)\,\widehat u(\mathbf k)
\bigr),
\]
computed using FFTs on the discrete wavenumber grid.
This single operator expression is used both for the physical source and for the
contrast--source terms arising during iteration.

Before the iteration begins, we form the incident field
$\psi_{\mathrm{inc}} = G_\varepsilon S$,
and construct the \emph{modified incident field}
\[
\psi_{\mathrm{mod}}^{(0)}(\mathbf r)
=
\frac{i}{\varepsilon}\,V(\mathbf r)\,G_\varepsilon S(\mathbf r),
\]
which serves as the initial iterate $\psi^{(0)}$.
This step requires only one application of $G_\varepsilon$ and one pointwise
multiplication and thus contributes negligibly to the total computational cost.

The CBS iteration is then carried out in real space, with the resolvent invoked
spectrally only for the nonlocal term.  
Using the operational form of the fixed-point map,
\[
\psi^{(n+1)}
=
\psi^{(n)}
+
\frac{i}{\varepsilon} V
\Bigl[
\psi_{\mathrm{inc}}
-
\psi^{(n)}
+
\mathcal F^{-1}\!\bigl(
\widehat g_\varepsilon(\mathbf k)\,
\widehat{\,V\psi^{(n)}\,}(\mathbf k)
\bigr)
\Bigr],
\]
the algorithm consists only of pointwise multiplications in real space and a
Fourier-accelerated convolution for $G_\varepsilon(V\psi^{(n)})$.
This hybrid real/spectral implementation mirrors the operator-theoretic analysis
of Section~4 and accommodates the complex-valued wavenumber $k(\mathbf r)^2$,
including the ABL, without modifying the background operator $H$. Apart from the presence of the ABL, our FFT-accelerated CBS iteration is similar to that of Osnabrugge et al. \cite{osnabrugge2016} as well as previous works on the CBS by Jakobsen et al. \cite{jakobsen2020convergent, jakobsen2020homotopy, jakobsen2020renormalized}.

All Fourier transforms are performed using FFTs on a uniform Cartesian grid.
For each frequency, the spectral kernel $\widehat g_\varepsilon(\mathbf k)$ is
precomputed once and stored in cache memory, so that each application of
$G_\varepsilon$ requires only a forward FFT, elementwise multiplication, and an
inverse FFT. Zero-padding is used to suppress wrap-around artefacts, after which
the solution is restricted to the physical domain.
Since the resolvent is implemented entirely as a spectral multiplier on a finite
wavenumber grid, no real-space regularization of $g_\varepsilon$ at the origin is
required.
 
The method is efficient for large-scale simulations: each iteration entails only
a few FFTs, giving $O(N\log N)$ complexity and excellent performance on modern
architectures \cite{jakobsen2020renormalized}.
Memory usage is low, as the algorithm stores only the current iterate, the
potential $V$, and the spectral kernel $\widehat g_\varepsilon$.  
The approach is also \emph{embarrassingly parallel} across frequencies: each
frequency uses its own kernel but shares FFT plans and data layout, enabling
highly scalable broadband simulations on CPUs, GPUs, and distributed-memory
systems.

\section{Numerical Results}

This section presents numerical experiments designed to assess the convergence
properties and practical performance of the operator-theoretic convergent Born
series (CBS) method. All experiments use the matrix-free FFT implementation
described in Section~6 together with the complex-density absorbing boundary layer
(ABL) introduced in Section~5.

\subsection{Homogeneous Medium}

We first consider a homogeneous medium with wave speed $c = 2200$~m/s. In the CBS
formulation, a homogeneous, non-dissipative reference medium with velocity
$c_0 = 2000$~m/s is used. Figure~1 shows the convergence history of the CBS
iterations. The method converges to a relative residual of $10^{-6}$ after
approximately 100 iterations.

Figure~2 compares frequency-domain wavefields computed using the CBS method, the
finite-difference frequency-domain (FDFD) method described in Appendix~C, and the
analytical Green’s function solution for a medium with velocity $c = 2200$~m/s.
Both the real and imaginary parts of the wavefields are visually indistinguishable
within the physical domain. Small differences are quantified in Figure~3, which
shows the pointwise amplitude and phase errors of the CBS and FDFD solutions
relative to the analytical Green’s function.

The amplitude and phase error patterns of the two numerical methods are very
similar. While the phase errors exhibit both positive and negative values, their
magnitudes remain small throughout the domain. Table~1 summarizes the overall
errors measured using a normalized $L^2$ norm. The CBS method yields a slightly
smaller overall error than the FDFD method, although both approaches achieve high
accuracy for this test case.

This result is noteworthy given that the CBS method employs a three-wavelength
viscoacoustic ABL, whereas the FDFD solution uses a standard perfectly matched
layer (PML) to suppress artificial boundary reflections. The close agreement with
the analytical solution demonstrates that the CBS formulation achieves accuracy
comparable to established frequency-domain solvers in homogeneous media.

\subsection{Heterogeneous medium with curved interfaces and strong contrasts}

We next assess the performance of the CBS method in a heterogeneous setting.
The test model consists of three layers separated by smoothly curved interfaces
(Figure~4), with P-wave velocities ranging from 2000 to 2800~m/s, representing
strong velocity contrasts.

Figure~5 shows the convergence behavior of the conventional Born series and the
CBS method for this model. While the conventional Born series diverges, the CBS
method converges robustly, reaching a relative residual of $10^{-6}$ after
approximately 200 iterations.

Frequency-domain wavefields computed using the CBS and FDFD methods are compared
in Figure~6. In both cases, absorption is handled using a three-wavelength
viscoacoustic absorbing buffer layer for CBS and a standard PML for FDFD. The
resulting wavefields are visually indistinguishable within the physical domain.
A more detailed comparison, shown in Figure~7, reveals small but measurable
differences in amplitude and phase between the two methods.

The overall discrepancy between the CBS and FDFD solutions, measured using a
normalized $L^2$ norm over the physical domain, is 2.13\%. This level of
agreement provides a numerical validation of the CBS method, demonstrating that
it reproduces the solution obtained with the standard FDFD approach at practical
accuracy levels.

Although FDFD is not free of the challenges discussed in the introduction, the
present experiment does not fall into the extreme high-frequency or strongly
scattering regime where those difficulties are most pronounced. Since FDFD
remains the standard forward solver in frequency-domain seismic full-waveform
inversion under the acoustic approximation, the close agreement observed here
supports the relevance of CBS for practical applications.

\section{Discussion}
\label{sec:discussion}

\subsection{Renormalized Forward Scattering Series}

A major motivation for this work is the use of scattering-series methods in
\emph{nonlinear inverse problems}. The inverse scattering series of Weglein
et al.~\cite{weglein2003inverse} provides a formally exact reconstruction
method, while Kouri and Vijay~\cite{kouri2003renormalization} argue that a renormalized Lippmann--Schwinger may be required for the development of an \emph{absolutely
convergent} inverse series. This renormalized framework has since been
extended to acoustic~\cite{lesage2013volterra} and multidimensional
settings~\cite{hussain2014volterra}. These developments emphasize that
convergence of the \emph{forward} scattering series is a necessary foundation
for any stable inverse construction, motivating the forward-scattering focus
of this subsection.

Although any linear system derived from the Lippmann--Schwinger equation can, in
principle, be solved using algebraic iterative methods such as GMRES, the existence of
a \emph{convergent scattering series} remains conceptually significant. Krylov methods
provide algebraic convergence but offer no physical interpretation of their iterates,
and their performance is highly sensitive to non-normality of the discretized operator.
Preliminary numerical experiments also suggest that embedding the CBS preconditioner
$\gamma$ within GMRES does not noticeably improve its convergence, indicating that the
CBS iteration is not primarily a preconditioning device but a \emph{renormalized,
physically interpretable forward scattering expansion}.

In classical scattering theory, the Born series represents the wavefield as a hierarchy
of single-, double-, and higher-order scattering events. This interpretation is simple
and transparent but notoriously fragile. As emphasized by Weinberg in his analysis of
the breakdown of the Born series in the presence of composite particles and resonances
\cite{Weinberg1963}, the series diverges whenever the underlying interaction produces
effective internal structure. In such situations, the bare interaction must be replaced
by a renormalized or ``dressed'' one. The CBS iteration embodies precisely this
principle: its powers $M^{n}$ serve as effective composite-scattering operators,
producing a renormalized forward series that remains convergent even in strongly
scattering media. This view connects CBS to renormalized Born and homotopy-based
constructions in acoustics and seismology
\cite{jakobsen2016renormalized,jakobsen2020renormalized}, while providing a rigorous
operator-theoretic guarantee of convergence.

A complementary interpretive perspective is obtained by viewing the CBS
damping parameter $\varepsilon$ as a scale parameter.
For clarity, this renormalization-group interpretation is formulated
strictly in the \emph{scalar} case, where the contrast acts as a
complex-valued multiplication operator and may be treated pointwise as a
scalar coupling.
In the operator-valued setting considered elsewhere in this paper, the
following discussion should therefore be regarded as heuristic.
Writing the shifted contrast as
\[
\Delta(\mathbf{r}) = V(\mathbf{r}) + i\varepsilon,
\]
we introduce the dimensionless effective coupling
\[
\bar{\lambda}(\varepsilon;\mathbf{r}) := \frac{\Delta(\mathbf{r})}{\varepsilon}.
\]
In the scalar setting, $\bar{\lambda}$ satisfies the simple flow relation
\[
\varepsilon \frac{d\bar{\lambda}}{d\varepsilon} = -\bar{\lambda},
\]
which makes explicit how increasing $\varepsilon$ suppresses strong local
contrast and enforces contractivity of the CBS iteration.
Using the Cayley-transform representation of the shifted resolvent, the
first two terms of the CBS expansion can be written entirely in terms of
$\bar{\lambda}$, making the resulting scale dependence explicit.
In this sense, increasing $\varepsilon$ suppresses fine-scale multiple
scattering interactions and enforces contractivity of the iteration.
The resulting behavior is formally analogous to the use of running
couplings in perturbative quantum field theory to control divergent
scattering expansions
\cite{GellMannLow1954,MandlShawQFT,Wilson1971}, although exploring such
analogies lies beyond the scope of the present study.
For our purposes, it is sufficient to view $\varepsilon$ as a scale
parameter and $\bar{\lambda}$ as a scale-dependent effective coupling.
The interpretation of $\varepsilon$ as a scale parameter is consistent
with the analysis of Osnabrugge et al., where $\varepsilon$ determines the
effective range of the dissipative Green’s function $G_\varepsilon$.
Likewise, interpreting $\bar{\lambda}$ as a scale-dependent (running)
coupling constant aligns naturally with the CBS convergence condition
$|\bar{\lambda}| < 1$.

\subsection{Extensions to Other Wave, Diffusion, and Open-System Problems}
\label{sec:extensions}

Although developed here for the scalar Helmholtz equation, the CBS convergence
analysis applies to a wider class of wave, diffusion, and open-system problems.
The essential requirements are that the background operator admit a stable
Cayley transform and that heterogeneity and dissipation enter through a bounded
multiplicative operator, ensuring a basis-independent numerical-range estimate.

Appendix~B illustrates this structure for the frequency-domain variable-density
acoustic equation \cite{XiangJakobsenGP2022}, where a first-order formulation
confines all derivatives to a self-adjoint background operator and places
parameter variations entirely in a bounded contrast. Similar augmented-state
reformulations exist for elastodynamics and related wave systems
\cite{jakobsenTransitionOperatorFWI}.

Dissipation enters naturally through the contrast operator, motivating the
complex-density absorbing boundary layer of Section~4.1, which preserves the
background operator and integrates directly with the convergence analysis.
The framework also extends to dissipative (non-self-adjoint) backgrounds that
admit a contractive Cayley transform, suggesting applicability to
low-frequency diffusion-dominated electromagnetic scattering problems
in conductive media 
\cite{saputeraEMdomainDecomposition}
as well as high-frequency electromagnetic scattering in dielectric media \cite{Krueger2017}. Since the Schroedinger equation can be expressed as a Helmholtz equation, it is in principle also possible to apply the CBS to open quantum systems \cite{Rotter2009NonHermitianOQS}.
Comprehensive numerical validation of these extensions is left for future work.

\section*{Concluding remarks}

This paper reformulates the Convergent Born Series (CBS) within a basis-independent operator-theoretic framework that yields a general convergence guarantee on bounded domains. By expressing the preconditioned Lippmann--Schwinger iteration solely in terms of the resolvent of a self-adjoint reference operator and its Cayley transform, the analysis avoids reliance on Fourier representations, translation invariance, or unbounded geometries, thereby removing key restrictions of the original CBS convergence theory.

A central consequence of this formulation is that dissipation and absorbing boundary layers can be incorporated directly into the convergence analysis by treating attenuation as part of a bounded multiplicative contrast operator. This separation preserves the self-adjoint structure of the background operator while strengthening contractivity of the CBS iteration, enabling the use of smoothly tapered complex-wavenumber absorbing layers within integral-equation solvers without modifying the differential operator. Numerical experiments demonstrate that solutions obtained in this manner closely agree with reference solutions computed using PML-based finite-difference frequency-domain methods.

From a computational perspective, the operator-theoretic formulation clarifies why the CBS iteration exhibits robust and predictable convergence in extended and strongly heterogeneous media, governed by an explicit local admissibility condition on the convergence parameter. Beyond forward modeling, the availability of a provably convergent forward scattering series provides a natural foundation for renormalized inverse scattering constructions. Although this work has focused on the scalar Helmholtz equation, the convergence analysis applies directly to a broader class of wave, diffusion, and open-system problems that admit a self-adjoint (or contractive) background operator and a bounded multiplicative contrast.

\section*{Declaration of generative AI and AI-assisted technologies}

During the preparation of this manuscript, the author used Microsoft Copilot as a supplementary tool for obtaining general overviews of relevant literature, for limited exploratory discussions, and for language editing and clarity improvements of text written by the author. GitHub Copilot was used to assist with streamlining the finite-difference frequency-domain (FDFD) code used for numerical benchmarking. All mathematical derivations, proofs, and final formulations were developed and independently verified by the author, who takes full responsibility for the content and correctness of the manuscript.

\clearpage 
    
\begin{table}[htbp]
    \centering
    \caption{Global relative $\ell^2$-errors in the physical domain for the homogeneous model.}
    \begin{tabular}{lcc}
        \hline \hline 
        Method & Relative $\ell^2$-error & Absorbing Layer \\
        \hline
        CBS   & 1.2141e-02 & ABL (constant-$Q$) \\
        FDFD  & 1.9814e-02 & PML (complex stretching) \\
        \hline \hline 
    \end{tabular}
    \label{tab:errors_hom}
\end{table}

\vspace{4cm}

\begin{figure}[h]
    \centering
    \includegraphics[width=0.9\textwidth]{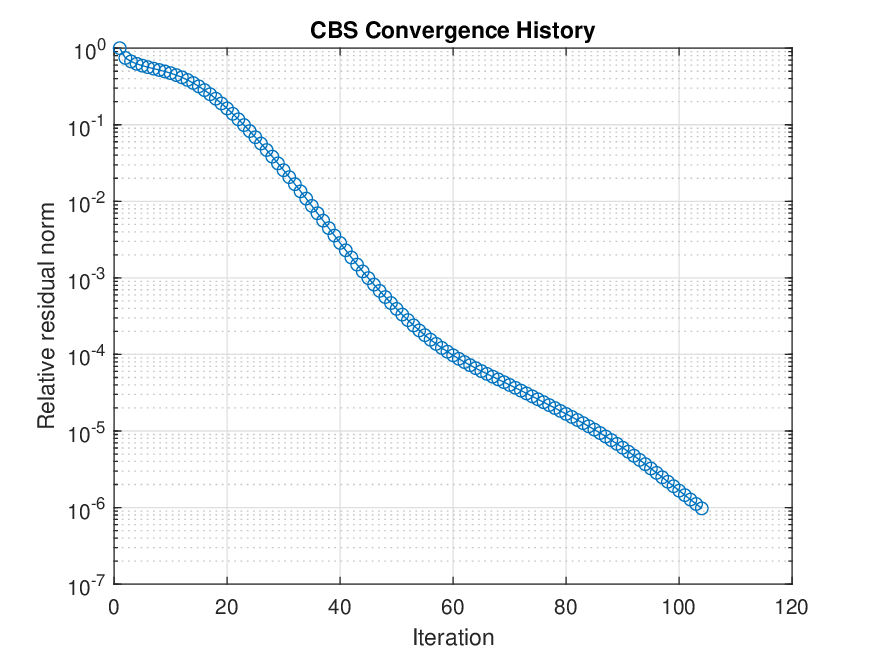}
    \caption{Convergence history of the CBS iteration in the homogeneous medium.}
    \label{fig:conv_hom}
\end{figure}

\begin{figure}[htbp]
    \centering
    \includegraphics[width=0.9\textwidth]{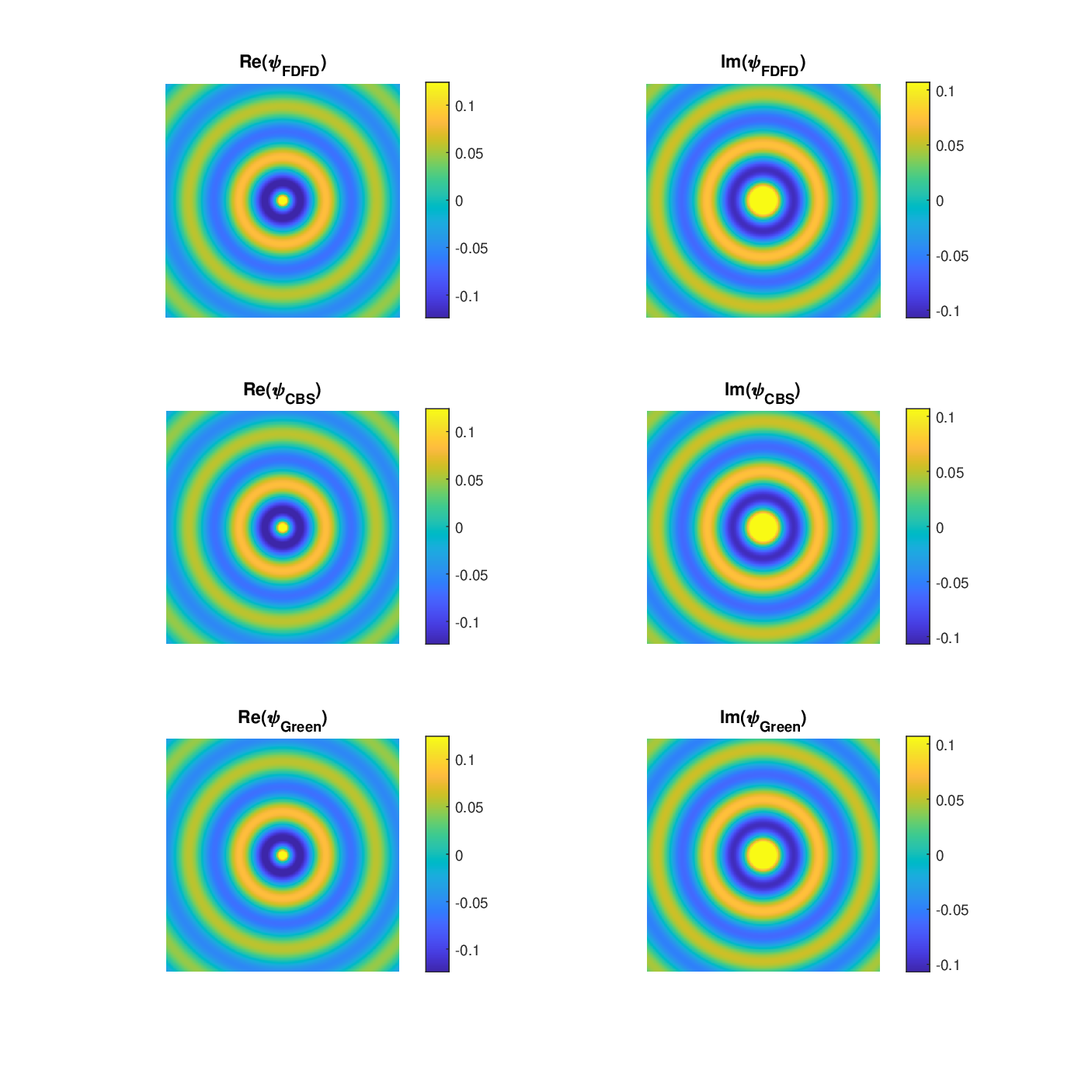}
    \caption{Real (left column) and imaginary (right column) parts of the wavefield within the physical domain of a homogeneous medium. Top: FDFD. Middle: CBS. Bottom: Green.}
    \label{fig:real_hom}
\end{figure}

\begin{figure}[htbp]
    \centering
    \includegraphics[width=0.9\textwidth]{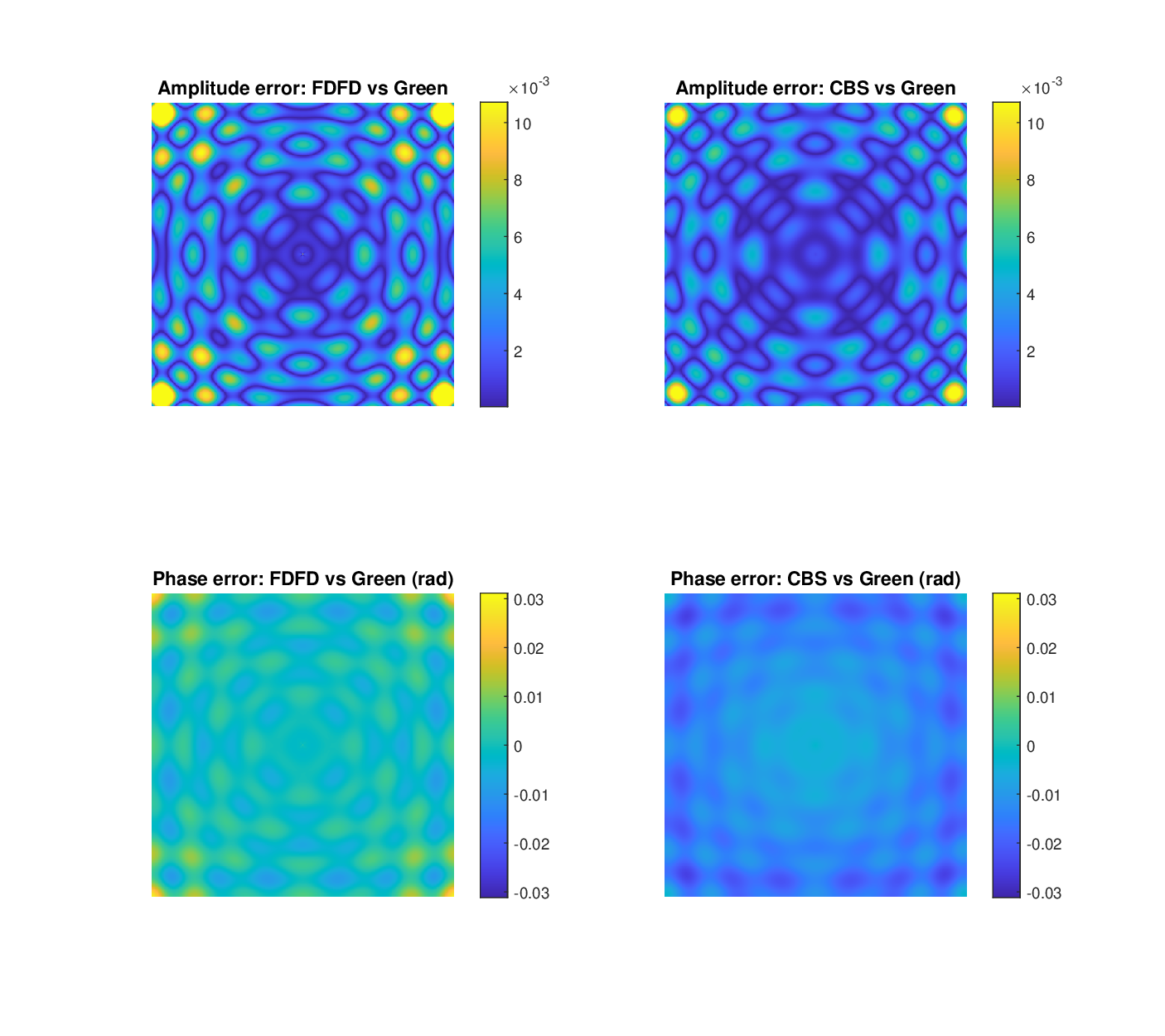}
    \caption{Relative amplitude (top) and absolute phase (bottom) errors of the FDFD (left) and CBS (right) relative to the analytical Green's function solution.}
    \label{fig:error_hom}
\end{figure}

\begin{figure}[htbp]
    \centering
    \includegraphics[width=\textwidth]{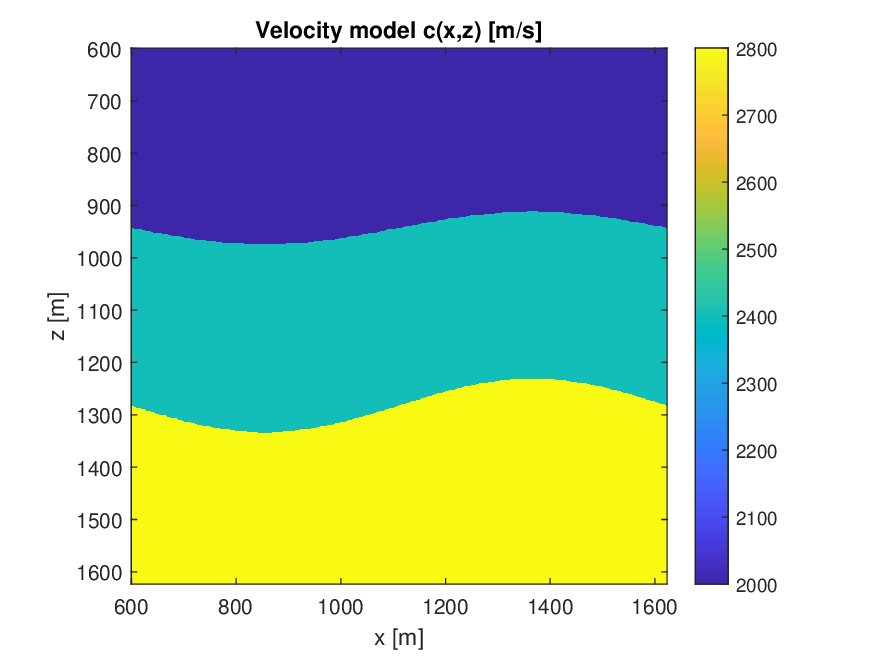}
    \caption{Heterogeneous velocity model used to compare CBS and FDFD results.}
    \label{fig:wavefields_imag_het}
\end{figure}

\clearpage 
\begin{figure}[htbp]
    \centering
    \includegraphics[width=0.9\textwidth]{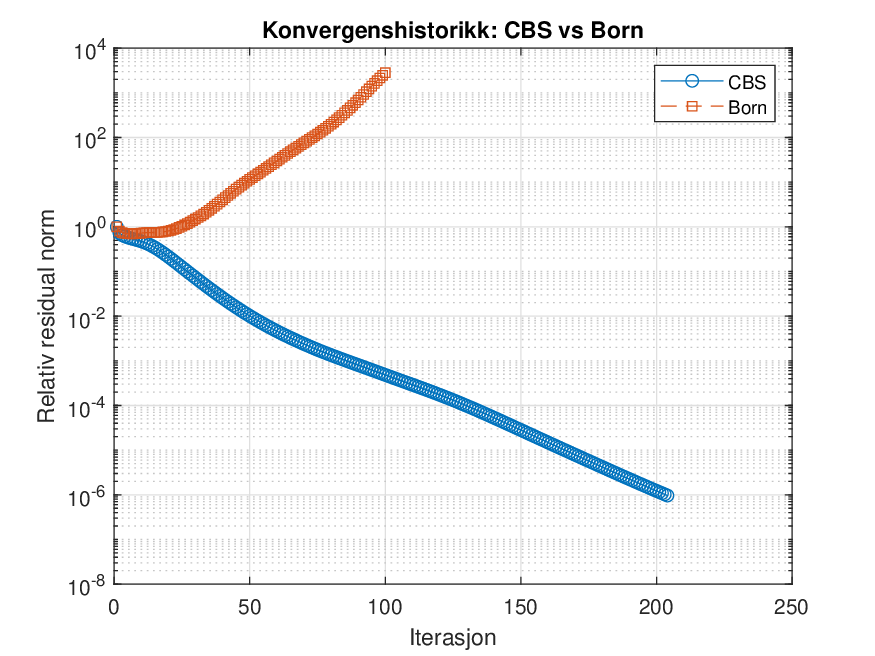}
    \caption{Convergence history of the CBS (blue curve) and the conventional Born series (red curve) for the heterogeneous medium shown in Figure 4. Note that the conventional Born series is diverging whereas the CBS is converging.}
    \label{fig:cbs_convergence_het}
\end{figure}

\clearpage 
\begin{figure}[htbp]
    \centering
    \includegraphics[width=\textwidth]{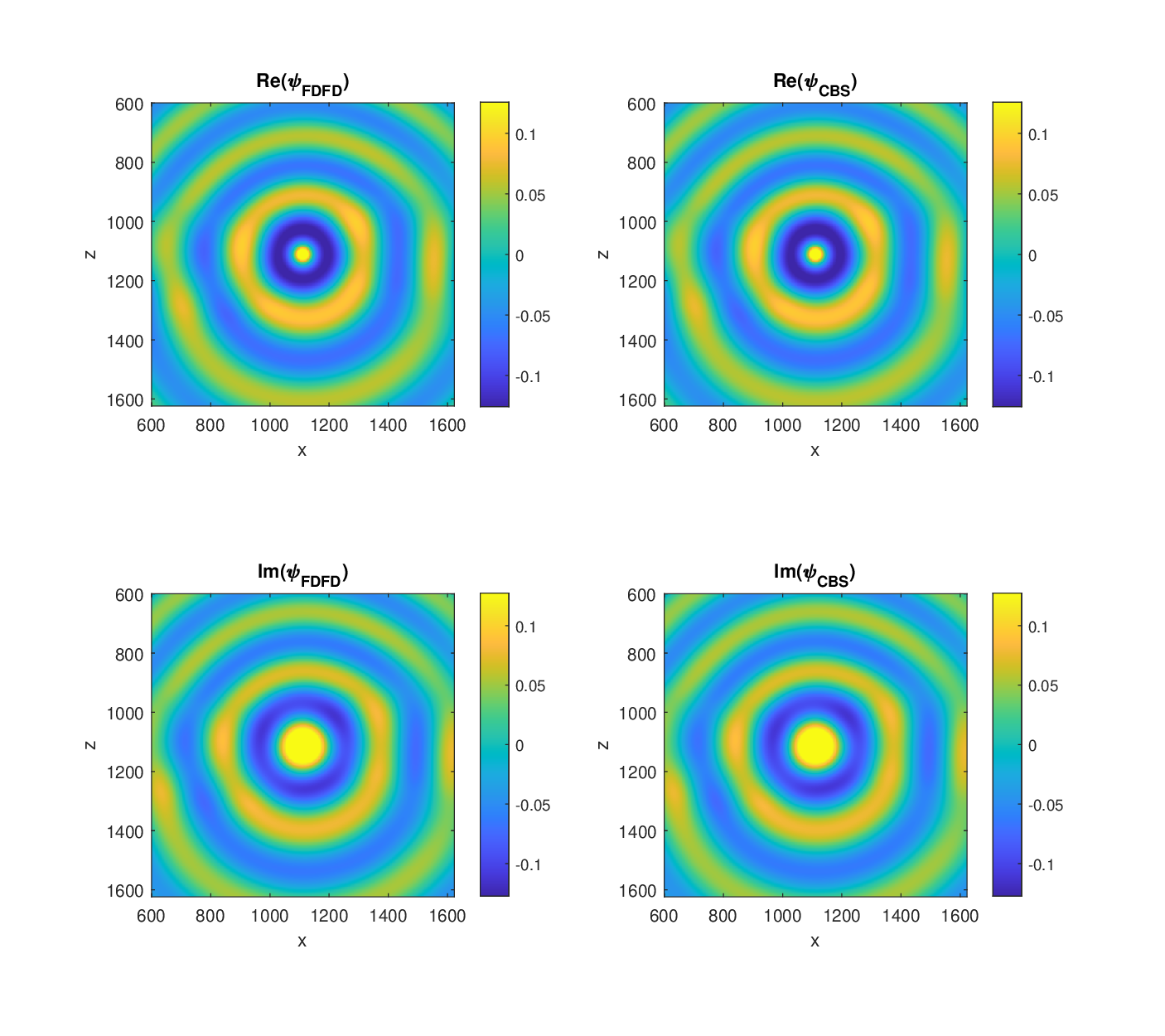}
    \caption{Frequency domain wavefields computed using the FDFD (left colum) and the CBS (red column). The real and imaginary parts of the wavefields are shown at the top and bottom, respectively.}
    \label{fig:relative_errors_het}
\end{figure}

\clearpage 
\begin{figure}[htbp]
    \centering
    \includegraphics[width=\textwidth]{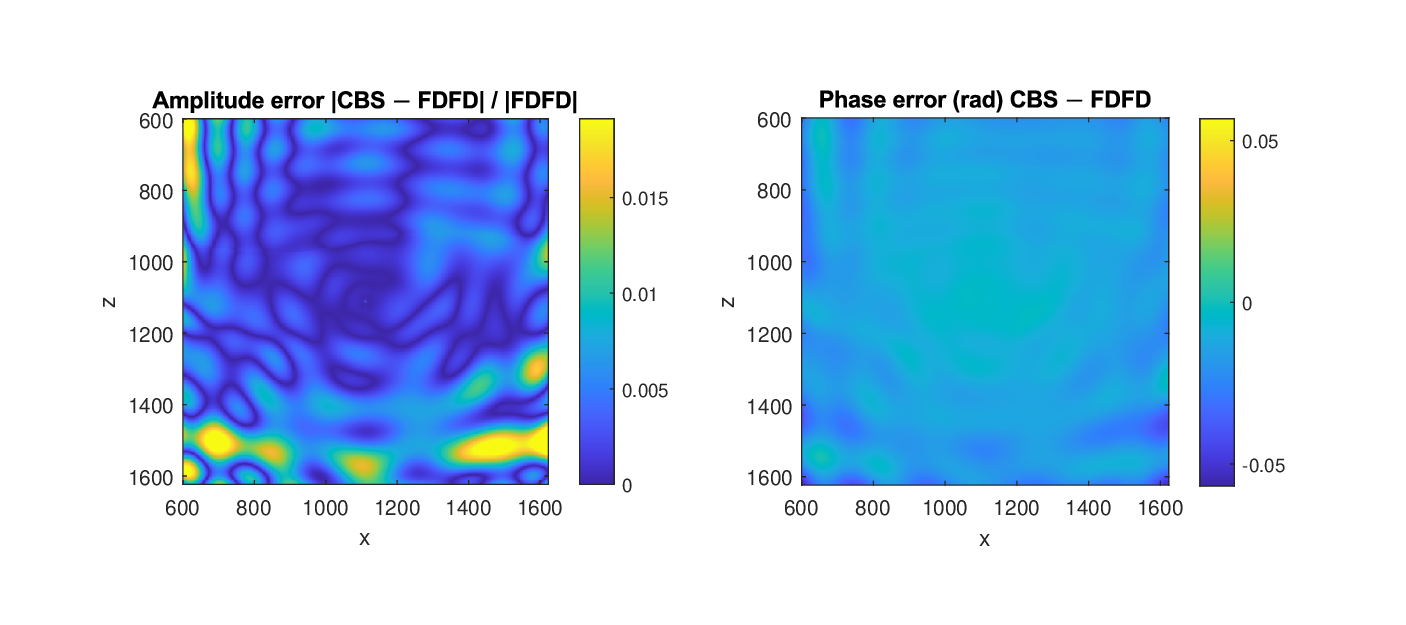}
    \caption{Amplitude difference between the wavefields computed using the CBS and FDFD methods. This figure represents a numerical validation of the CBS method using the more common FDFD method as the reference.}
    \label{fig:relative_errors_het}
\end{figure}

\clearpage 
\bibliographystyle{elsarticle-num}
\bibliography{references}


\appendix

\clearpage 

\section{Finite Difference Frequency Domain Method with PML}

To benchmark the performance of the CBS method, we implemented a finite difference frequency domain (FDFD) solver for the Helmholtz equation using an advanced perfectly matched layer (PML) formulation. This appendix summarizes the numerical scheme, boundary treatment, and validation strategy. Relevant references are \cite{marfurt1984, ajo2005frequency, AbubakarHabashy2013ViscoAcousticVIE}.

\subsection{FDFD Formulation and Linear System}

We discretize the scalar Helmholtz equation
\[
\nabla^2 \psi(\mathbf{r}) + k^2(\mathbf{r}) \psi(\mathbf{r}) = -S(\mathbf{r})
\]
on a regular Cartesian grid using a 9-point stencil for the Laplacian operator. This leads to a sparse linear system of the form
\[
A \psi = b,
\]
where $A$ is a complex-valued sparse matrix representing the discretized differential operator and the spatially varying wavenumber $k(\mathbf{r})$, $\psi$ is the unknown wavefield vector, and $b$ is the source term vector. We refer to the unpublished paper by Franklin \cite{ajo2005frequency} for a good discussion of the FDFD method for solving the Helmholtz equation. The paper of Huang and Greenhalgh \cite{Huang} also discuss relevant details and is particularly interesting here since it present a kind of FDFD version of the CBS method. 

In our implementation, this system is solved in MATLAB using the built-in
backslash operator,
\[
\texttt{psi = A \textbackslash\ b}.
\]
For sparse, square, complex-valued matrices such as those arising from
FDFD discretizations of the Helmholtz equation with PML, MATLAB automatically
applies a direct sparse LU factorization with partial pivoting. Internally,
this involves fill-reducing reorderings and robust numerical pivoting
strategies, providing a stable and reproducible solution without the need for
user-specified solver parameters.

Although direct solvers are not scalable to very large problem sizes, they
are well suited for reference computations on moderate grids and therefore
provide a reliable baseline for benchmarking the CBS method.

\subsection{Perfectly Matched Layer}

To suppress artificial reflections from domain boundaries, we employ a complex coordinate stretching approach \cite{berenger1994pml, chew1996pml}, where the spatial derivatives are modified as
\[
\frac{\partial}{\partial x} \rightarrow \frac{1}{1 + i\sigma_x(x)/\omega} \frac{\partial}{\partial x},
\]
with $\sigma_x(x)$ denoting the damping profile in the PML region and $\omega$ the angular frequency. The damping profile is chosen to vary smoothly from zero at the interface to a maximum value at the outer boundary, ensuring minimal reflection for all angles of incidence.

\clearpage 
\section{Convergent Born Series for Variable-Density Acoustics}
\label{app:AcousticCBS}

This appendix is included solely to demonstrate that the operator-theoretic assumptions underlying the Cayley-transform convergence analysis developed in Sections~3--4 are not specific to the scalar constant-density Helmholtz equation. In particular, we show that the frequency-domain variable-density acoustic wave equation can be reformulated in a first-order form in which all spatial derivatives are confined to a self-adjoint background operator, while heterogeneity and dissipation enter exclusively through a bounded, purely multiplicative contrast operator. This establishes structural admissibility of variable-density acoustics within the CBS framework at the level of operator theory. No claims are made here regarding numerical performance, discretization choices, or computational efficiency, and no additional numerical experiments are presented. This appendix is only included to illustrate the generality of the convergence argument and to clarify the class of wave equations to which the Cayley-transform analysis applies in principle.

\subsection{First-order formulation}

We consider time-harmonic acoustic wave propagation with time dependence
$e^{-i\omega t}$ in a heterogeneous fluid characterized by spatially varying
density $\rho(\mathbf r)$ and compressibility $\kappa(\mathbf r)=1/K(\mathbf r)$.
The frequency-domain acoustic equations may be written in first-order form as \cite{deHoop1995Handbook, XiangJakobsenGP2022, JakobsenXiangVanDongen2023, JakobsenSaputeraRadu2025}
\begin{equation}
\begin{aligned}
-i\omega\,\rho(\mathbf r)\,\mathbf v(\mathbf r) &= \nabla p(\mathbf r), \\
-i\omega\,\kappa(\mathbf r)\,p(\mathbf r) &= \nabla\cdot\mathbf v(\mathbf r)
+ s(\mathbf r),
\end{aligned}
\label{eq:first_order_acoustics}
\end{equation}
where $p$ denotes pressure, $\mathbf v$ particle velocity, and $s$ an impressed
source.

Introducing the state vector
\[
\psi(\mathbf r)
=
\begin{pmatrix}
\mathbf v(\mathbf r) \\
p(\mathbf r)
\end{pmatrix}
\in L^2(\Omega;\mathbb C^4),
\]
equation~\eqref{eq:first_order_acoustics} can be written compactly as
\begin{equation}
(A + V(\mathbf r))\,\psi = f,
\label{eq:acoustic_operator}
\end{equation}
with $f=(0,s)^T$.

\subsection{Self-adjoint reference operator}

We introduce constant reference parameters $\rho_0>0$ and $\kappa_0>0$ and define
the reference operator
\begin{equation}
A :=
\begin{pmatrix}
0 & \nabla \\
\nabla\cdot & 0
\end{pmatrix}.
\label{eq:A_acoustics}
\end{equation}
After multiplication of~\eqref{eq:acoustic_operator} by $i$, the reference
operator $A$ is self-adjoint on $L^{2}(\Omega;\mathbb{C}^{4})$ under standard
boundary conditions, and all spatial derivatives appear exclusively in $A$;
equivalently, one may rescale the velocity and pressure variables by constant
reference parameters so that $A$ is self-adjoint with respect to the standard
$L^{2}$ inner product, the two formulations being unitarily equivalent.

\subsection{Multiplicative contrast operator}

Writing
\[
\rho(\mathbf r) = \rho_0 + \delta\rho(\mathbf r),
\qquad
\kappa(\mathbf r) = \kappa_0 + \delta\kappa(\mathbf r),
\]
the contrast operator takes the diagonal, purely multiplicative form
\begin{equation}
V(\mathbf r) =
\begin{pmatrix}
-i\omega\,\delta\rho(\mathbf r)\,I_3 & 0 \\
0 & -i\omega\,\delta\kappa(\mathbf r)
\end{pmatrix}.
\label{eq:V_acoustics}
\end{equation}
The operator $V$ is bounded on $L^2(\Omega;\mathbb C^4)$ and contains no
derivatives, which is the key structural requirement for the CBS convergence
theory.

\subsection{Shifted resolvent and admissibility condition}

Following Sections~3--4, we introduce the shifted resolvent
\[
G_\varepsilon := (A + i\varepsilon I)^{-1},
\qquad \varepsilon > 0,
\]
and the associated Cayley transform
\[
U := (A - i\varepsilon I)(A + i\varepsilon I)^{-1}.
\]
Because $A$ is self-adjoint, $U$ is unitary. The CBS convergence analysis of
Section~3.3 applies directly, and geometric convergence of the CBS iteration is
guaranteed provided the admissibility condition
\begin{equation}
\varepsilon \ge \|V\|_{L^\infty}
\label{eq:acoustic_admissibility}
\end{equation}
is satisfied.

\subsection{Lippmann--Schwinger equation and CBS iteration}

Applying $G_\varepsilon$ to~\eqref{eq:acoustic_operator} yields the
Lippmann--Schwinger equation
\begin{equation}
\psi = G_\varepsilon f + G_\varepsilon V \psi,
\label{eq:LS_acoustics}
\end{equation}
\emph{without a minus sign}. This sign convention is consistent with
Sections~2--4 and reflects the decomposition $A+V$ used in the CBS framework.

Using the standard CBS preconditioning, the fixed-point iteration takes the form
\begin{equation}
\psi^{(n+1)} = M\psi^{(n)} + b,
\end{equation}
with
\begin{equation}
M = I - \frac{i}{\varepsilon}V
    + \frac{i}{\varepsilon}V G_\varepsilon V,
\qquad
b = \frac{i}{\varepsilon}V G_\varepsilon f.
\label{eq:CBS_acoustics}
\end{equation}
Under the admissibility condition~\eqref{eq:acoustic_admissibility}, the numerical
range of $M$ lies within the closed unit disk, and the CBS iteration converges
geometrically.

\end{document}